\documentclass[12pt,amsfonts]{article}

\usepackage{latexsym}
\usepackage{amsmath}
\usepackage{amsfonts,amssymb} 

\newtheorem{theorem}{Theorem}[section]
\newtheorem{lemma}[theorem]{Lemma}%[section]
\newtheorem{corollary}[theorem]{Corollary}%[section]
\newtheorem{proposition}[theorem]{Proposition}%[section]

\newtheorem{definition}[theorem]{Definition}
\newcommand{\bd}[1]{\begin{definition}\label{#1}\rm}
\newcommand{\ed}{\end{definition}}
\newcommand{\bt}[1]{\begin{theorem}\label{#1}}
\newcommand{\et}{\end{theorem}}
\newcommand{\bprop}[1]{\begin{proposition}\label{#1}}
\newcommand{\eprop}{\end{proposition}}
\newcommand{\bcor}[1]{\begin{corollary}\label{#1}}
\newcommand{\ecor}{\end{corollary}}

\newcommand{\D}{\displaystyle}

\newcommand{\vp}{\varphi}
\newcommand{\ve}{\varepsilon}
\newcommand{\nid}{\noindent}
\newcommand{\qed}{\hfill$\Box$} 

\def\1{\, {\rm I}\mskip-10mu 1} 
\,

\renewcommand{\t}[1]{\tilde{#1}} 
\newcommand{\bmu}{\mbox{\boldmath${\mu}$}}
\newcommand{\bnu}{\mbox{\boldmath${\nu}$}} 
\newcommand{\sbnu}{\mbox{\tiny\boldmath${\nu}$}} 
\newcommand{\sbmu}{\mbox{\tiny\boldmath${\mu}$}} 

\begin{document}
\title{{Absolute Continuity under Time Shift for Ornstein-Uhlenbeck type Processes with 
Delay or Anticipation}}  
\par
\author{J\"org-Uwe L\"obus
%%\footnotemark[1]
\\ Matematiska institutionen \\ 
Link\"opings universitet \\ 
SE-581 83 Link\"oping \\ 
Sverige 
}
\date{}
\maketitle
{\footnotesize
\noindent
\begin{quote}
{\bf Abstract}
%\footnotetext[1]{\noindent } 
The paper is concerned with one-dimensional two-sided Ornstein-Uhlenbeck 
type processes with delay or anticipation. We prove existence and uniqueness 
requiring almost sure boundedness on the left half-axis in case of delay 
and almost sure boundedness on the right half-axis in case of anticipation. 
For those stochastic processes $(X,P_{\sbmu})$ we calculate the 
Radon-Nikodym density under time shift of trajectories, $P_{\sbmu}(dX_{\cdot 
-t})/P_{\sbmu}(dX)$, $t\in {\Bbb R}$. 
\noindent 

{\bf AMS subject classification (2010)} primary 60J65, secondary 60H10

\noindent
{\bf Keywords} 
Brownian motion, SDE with delay or anticipation, non-linear transformation 
of measures 
\end{quote}
}
%\bigskip
%\noindent (Running title: )

\section{Introduction, Basic Objects, and Main Result}
\setcounter{equation}{0}

\subsection{Introduction}

Let $W$ be a one-dimensional two-sided Brownian motion with random 
initial value $W_0$ and $a\in (-1,0)$ be non-random. The paper has two 
objectives. On the one hand, we are interested in solutions to a 
stochastic equation with delay of the form  
\begin{eqnarray*}
dX_s=a\left(X_{s-1}+b_0\vphantom{l^1}\right)\, ds+dW_s\, ,\quad s\in {\Bbb 
R},\quad {\rm where}\ \limsup_{v\to -\infty}\left|X_v\right|<\infty, 
\end{eqnarray*} 
and solutions to a stochastic equation with anticipation of the form 
\begin{eqnarray*}
dX_s=-a\left(X_{s+1}+b_0\vphantom{l^1}\right)\, ds+dW_s\, ,\quad s\in {\Bbb 
R},\quad{\rm where}\ \limsup_{v\to +\infty}\left|X_v\right|<\infty, 
\end{eqnarray*} 
for some $b_0\equiv b_0(W)$, cf. Theorem \ref{Theorem3.1} and Corollary 
\ref{Corollary3.2} below. This part of the paper can be regarded as a 
supplement to the rich literature on stochastic delay differential 
equations. We wish to point to the many different topics of interest by 
referring to some of them, \cite{GK00}, \cite{MS90}, \cite{NR07}, and 
\cite{Ri08}. 

On the other hand, for the stochastic processes $(X,P_{\sbmu})$ being solution 
to one of the two above equations, we are interested in quasi-invariance under 
time shift of trajectories. In Theorem \ref{Theorem3.3} we provide the 
Radon-Nikodym density $P_{\sbmu}(dX_{\cdot -t})/P_{\sbmu}(dX)$, $t\in {\Bbb R}$, 
and relate this density to a recent result in the paper \cite{Lo13}. This type 
of quasi-invariance is particularly meaningful for two-sided stochastic 
processes. It adds to the well-established results on quasi-invariance in the 
non-anticipating or the non-gaussian case. For an overview on these, check for 
example, \cite{B10}, Section 11.2, \cite{Bu94}, and \cite{UZ00}. 

\subsection{Definitions and Notations} 

{\bf Two-sided Brownian motion with random initial value. } 
Let $\bnu$ be a probability distribution on $({\Bbb R},{\mathcal B}({\Bbb R}))$. 
Let the probability space $(\Omega,{\mathcal F},Q_{\sbnu})$ be given by the 
following. For a similar definition see also \cite{Im96}, paragraph 1 of Section 2.  
\begin{itemize} 
\item[{(i)}] $\Omega =C({\Bbb R};{\Bbb R})$, the space of all continuous 
functions $W$ from ${\Bbb R}$ to ${\Bbb R}$. Identify $C({\Bbb R};{\Bbb R}) 
\equiv\{\omega\in C({\Bbb R};{\Bbb R}):\omega(0)=0\}\times {\Bbb R}$. 
%\item[{}] \vspace{-1.1cm}
\item[{(ii)}] ${\mathcal F}$, the $\sigma$-algebra of Borel sets with respect 
to uniform convergence on compact subsets of ${\Bbb R}$.  
%\item[{}] \vspace{-1.1cm}
\item[{(iii)}] $Q_{\sbnu}$, the probability measure on $(\Omega ,{\mathcal F})$ 
for which, when given $W_0$, both $(W_s-W_0)_{s\ge 0}$ as well as $(W_{-s}-
W_0)_{s\ge 0}$ are independent standard Brownian motions with state space 
$({\Bbb R},{\mathcal B}({\Bbb R}))$. Assume $W_0$ to be independent of $(W_s-W_0 
)_{s\in {\Bbb R}}$ and distributed according to $\bnu$. 
\end{itemize} 
In addition, we will assume that the natural filtration $\{{\mathcal F}_u^v= 
\sigma(W_\alpha-W_\beta :u\le\alpha,\beta\le v)\times \sigma(W_0):-\infty 
<u<v<\infty\}$ is completed by the $Q_{\sbnu}$-completion of ${\mathcal F}$. 
\medskip 

For the measure $\bnu$ we shall assume the following throughout the paper.  
\begin{itemize}
\item[{(i)}] $\bnu$ admits a density $m$ with respect to the one-dimensional 
Lebesgue measure. 
\item[{(ii)}] $0<m\in C^1({\Bbb R};{\Bbb R})$.  
\end{itemize} 
It is known that 
\begin{eqnarray}\label{1.1}
\frac{Q_{\sbnu}(dW_{\, \cdot -t})}{Q_{\sbnu}(dW)}=\frac{m\left(W_{-t} 
\vphantom{l^1}\right)}{m(W_0)}\, , \quad t\in {\Bbb R}.  
\end{eqnarray} 
In \cite{Lo13}, Theorem 1.12 together with Proposition 1.16, this formula 
is generalized. Below we will outline the result of this reference in a 
simplified form for which we assume for technical reasons 
\begin{itemize}
\item[{(iii)}] that for all $y\in {\Bbb R}$,  
\begin{eqnarray*} 
\left.\frac{d}{d\lambda}\right|_{\lambda =0}\frac{m(\, \cdot +\lambda y)} 
{m}=y\cdot\frac{\nabla m}{m}\quad \mbox{\rm exists in } L^Q({\Bbb R},\bnu) 
\end{eqnarray*} 
for some $Q\in (1,\infty)$. 
\end{itemize} 

\nid
{\bf Parallel trajectories. } Let $\xi\equiv\xi(W)$ be a random variable 
and let $\1\equiv\1(s)$, $s\in {\Bbb R}$, denote the constant function on 
${\Bbb R}$ taking the value one. For fixed $W\in\Omega$ and $x\in {\Bbb R}$, 
denote  
\begin{eqnarray*} 
\nabla_x\xi(W+x\1)\equiv\frac{\partial}{\partial x}\xi(W+x\1)
\end{eqnarray*} 
whenever this partial derivative exists. Set 
\begin{eqnarray*} 
\nabla_{W_0}\xi(W):=\left.\nabla_x\right|_{x=W_0}\xi(W-W_0\1+x\1)\, . 
\end{eqnarray*} 

\nid
{\bf The process $X$. } In Theorem \ref{Theorem3.1} and Corollary 
\ref{Corollary3.2} below, we will verify that the stochastic processes of 
Subsection 1.1 have a representation 
\begin{eqnarray*} 
X=W+A 
\end{eqnarray*} 
where $A\equiv A(W)$ is a stochastic process with trajectories that 
belong to $C({\Bbb R};{\Bbb R})$. This induces path wise a measurable map 
$X=X(W):\Omega\equiv C({\Bbb R};{\Bbb R})\to C({\Bbb R};{\Bbb R})$ which we 
show to be $Q_{\sbnu}$-a.e. injective. The measure $P_{\sbmu}:=Q_{\sbnu} 
\circ X^{-1}$ is the law of $X$. Here, the index $\bmu$ is just a part of a 
symbol added for compatibility of the notation with \cite{Lo13}. 

\subsection{Change of measure under time shift}

Following the ideas of \cite{Lo13}, we shall verify in Theorem 
\ref{Theorem3.1} and Corollary \ref{Corollary3.2} below that $X\equiv X(W)$ 
is a {\it temporally homogeneous} function of $W\in\Omega$. In the present 
setting this means that for $W^u:=W_{\cdot +u}+A_u(W)\1$, $u\in {\Bbb R}$, 
we have $A_0(W)=0$ and 
\begin{eqnarray*} 
X_{\cdot+v}(W)=X(W^v)\, , \quad v\in {\Bbb R}. 
\end{eqnarray*} 
In this case we have, among other things, $X_0=W_0$. We may therefore also 
write $\nabla_{X_0}$ for $\nabla_{W_0}$. 
\medskip

One purpose of the present paper is to show that the Radon-Nikodym densities  
\begin{eqnarray}\label{1.2} 
\omega_{-t}(W):=\frac{Q_{\sbnu}(dW^{-t})}{Q_{\sbnu}(dW)}\quad\mbox{\rm and} 
\quad\rho_{-t}(X):=\frac{P_{\sbmu}(dX_{\cdot -t})}{P_{\sbmu}(dX)}\, , \quad 
t\in {\Bbb R}, 
\end{eqnarray} 
exist and to give a representation of these densities. In fact, in Section 
3, we will use a particular property of the processes we are examining to 
establish the above Radon-Nikodym derivatives. However we will also discuss 
the relation to one the main results of \cite{Lo13}, namely Theorem 1.12. 
There, for a much larger class of processes than those of Subsection 1.1, 
we derive $Q_{\sbnu}$-a.e. the density 
\begin{eqnarray*} 
\omega_{-t}(W)=\frac{m\left(X_{-t}(W)\vphantom{l^1}\right)}{m(W_0)}\cdot 
\left|1+\nabla_{W_0}A_{-t}(W)\vphantom{\dot{f}}\right|  
\end{eqnarray*} 
and $P_{\sbmu}$-a.e. the density 
\begin{eqnarray}\label{1.3}
\rho_{-t}(X)=\frac{\D m\left(X_{-t}\vphantom{l^1}\right)}{\D m\left(X_0 
\vphantom{l^1}\right)}\cdot\vphantom{\left\{\sum_{\tau_k}\right\}}\left| 
\nabla_{X_0}X_{-t}\vphantom{\dot{f}}\right|\, ,\quad t\in {\Bbb R}. 
\end{eqnarray} 

\section{Ornstein-Uhlenbeck type process with delay or anticipation} 
\setcounter{equation}{0}

In order to approach solutions to the stochastic equations of Subsection 
1.1 we are first interested in a process $\t X$ that obeys a stochastic 
delay equation of the form 
\begin{eqnarray}\label{2.1}
d\t X_s=a\t X_{s-1}\, ds+dW_s\, ,\quad s\in {\Bbb R},\quad\t X_u=\vp(u,W) 
\, ,\quad u\in [-1,0], 
\end{eqnarray} 
or a stochastic equation with anticipation of the form 
\begin{eqnarray}\label{2.2}
d\t X_s=-a\t X_{s+1}\, ds+dW_s\, ,\quad s\in {\Bbb R},\quad\t X_u=\vp(u,W) 
\, ,\quad u\in [0,1],  
\end{eqnarray} 
for a non-random $a\in (-1,0)$ and a certain continuous random function 
$\vp$. Our approach will rely on a very specific choice of $\vp$, cf. 
Theorem \ref{Theorem2.5} and Corollary \ref{Corollary2.6} below. In fact 
the construction of such a {\it two-sided} process $\t X$ will suggest to 
relate uniqueness to the asymptotic behavior rather than to an initial 
function.  
\medskip 

Our first effort will be towards the $Q_{\sbnu}$-a.e. pathwise existence and 
uniqueness of the solution to the equation (\ref{2.1}). It will turn out that 
the $Q_{\sbnu}$-a.e. pathwise existence and uniqueness of the solution to 
(\ref{2.2}) is a consequence of the result for (\ref{2.1}), simply by time 
reversal. 
\medskip

Our method to treat the equation (\ref{2.1}) is adapted to the {\it method of 
steps} in the theory of delay differential equations. To begin with, we are going 
to specify the space of functions the initial function $\vp$ is related with. 
\begin{definition}\label{Definition2.1}
(a) For $a\in {\Bbb R}\setminus\{0\}$, $\omega\in C((-\infty,0])$, and $k\in 
{\Bbb N}$, let $C_{a,\omega}^k([-1,0])$ denote the space of all real functions 
$f$ on $[-1,0]$ satisfying the following. For all those $f$ there is a function 
$p_k\in C([-1,0])$ and a collection $c_0,c_1,\ldots ,c_k\in {\Bbb R}$ such that 
$f$ has a representation in the form of a $k$-fold iterated integral  
\begin{eqnarray*} 
f(v_0)&&\hspace{-.5cm}=\int_{-1}^{v_0} a\left(\ldots a\left(\int_{-1}^{v_{k-1}} 
a\left(p_k(v_{k})+c_k+\omega (v_k-k)\vphantom{\dot{f}}\right)\, dv_k+c_{k-1} 
\right.\right. \\ 
&&\hspace{-1.2cm}\left.\left.+\omega(v_{k-1}-(k-1))\vphantom{\int_0^{v_{k-1}}} 
\right)\, dv_{k-1}\ldots +c_1+\omega (v_1-1)\right)\, dv_1+c_0+\omega (v_0)\, , 
\end{eqnarray*} 
$v_0\in [-1,0]$. In addition, let $p_0\in C([-1,0])$ be defined by 
\begin{eqnarray*}
p_0:=f-(c_0+\omega )\, . 
\end{eqnarray*} 
(b) 
\begin{eqnarray*}
C_{a,\omega}^\infty([-1,0]):=\bigcap_{k\in {\Bbb N}}C_{a,\omega}^k([-1,0])\, . 
\end{eqnarray*} 
\end{definition} 
Obviously, 
\begin{eqnarray*}
C_{a,\omega}^{k+1}([-1,0])\subseteq C_{a,\omega}^k([-1,0])\, , \quad k\in  
{\Bbb N}. 
\end{eqnarray*} 
Furthermore, for $f\in C_{a,\omega}^k([-1,0])$ represented as in part (a) 
of this definition, let us use the notation $f=T^{(1)}(p_1;c_0,c_1)=T^{(2)} 
(p_2;c_0,c_1,c_2)=\ldots =T^{(k)}(p_k;c_0,\ldots$ $,c_k)$ and $p_i\equiv p_i 
(\cdot ,f)$, $i=1,\ldots ,k$, $c_j\equiv c_j(f)$, $j=0,\ldots ,i$. In the 
proof below we will also use the symbol $C_{a,\omega(\cdot -k)}^1([-1,0])$ 
for $\omega$ replaced by the restriction of $\omega(\cdot -k)$ to $(-\infty, 
0]$, $k\in {\Bbb N}$. 
\begin{proposition}\label{Proposition2.2} 
Let $a\in {\Bbb R}\setminus\{0\}$ and $\omega\in C((-\infty,0])$. (a) For 
all $k\in {\Bbb N}$, $C_{a,\omega}^k([-1,0])$ is dense in $C([-1,0])$. \\ 
(b) The set $C_{a,\omega}^\infty([-1,0])$ is dense in $C([-1,0])$. 
\end{proposition}
Proof. (a) This follows immediately by induction. 
\medskip 

\nid 
(b) Let $\vp\in C([-1,0])$, $0<\ve <\frac12$, and $\vp_1\in C_{a,\omega 
}^1([-1,0])$ such that $\|\vp -\vp_1\|<\ve$ where $\|\cdot\|$ denotes 
the norm in $C([-1,0])$ in this proof. Furthermore, choose $\psi_1$ with 
$\psi_1(-1)=0$ and $\psi_1+\omega(\cdot -1)|_{[-1,0]}\in C_{a,\omega( 
\cdot -1)}^1([-1,0])$ such that $\|\psi_1-p_1(\cdot ,\vp_1)\|<\ve^2\wedge 
\ve^2/|a|$. Setting $\vp_2:=T^{(1)}(\psi_1;c_0(\vp_1),c_1(\vp_1))$ it 
follows that 
\begin{eqnarray*}
\vp_2\in C_{a,\omega}^2([-1,0])
\end{eqnarray*} 
and 
\begin{eqnarray*}
\|p_1(\cdot ,\vp_1)-p_1(\cdot ,\vp_2)\|<\ve^2\quad\mbox{\rm as well as} 
\quad\|\vp_1-\vp_2\|<\ve^2\, . 
\end{eqnarray*} 
We continue that way by choosing for all $k\in {\Bbb N}$ a $\psi_k$ with 
$\psi_k(-1)=0$ and $\psi_k+\omega(\cdot -k)|_{[-1,0]}\in C_{a,\omega(\cdot 
-k)}^1([-1,0])$ such that $\|\psi_k-p_k(\cdot ,\vp_k)\|<\ve^{k+1}\wedge 
\ve^{k+1}/|a|^k$. We put $\vp_{k+1}:=T^{(k)}(\psi_k;c_0(\vp_k),\ldots ,c_k 
(\vp_k))$ and obtain 
\begin{eqnarray}\label{2.3} 
\vp_{k+1}\in C_{a,\omega}^{k+1}([-1,0])\quad\mbox{\rm and}\quad\|\vp_k- 
\vp_{k+1}\|<\ve^{k+1} 
\end{eqnarray} 
as well as 
\begin{eqnarray*}
\|p_i(\cdot ,\vp_k)-p_i(\cdot ,\vp_{k+1})\|<\ve^{k+1}\, , \quad i\in \{1, 
\ldots ,k\}\, . 
\end{eqnarray*} 
From here it follows that the sequences $(\vp_k)_{k\in {\Bbb N}}$ and $( 
p_i(\cdot ,\vp_k))_{k\ge i}$, $i\in {\Bbb N}$, simultaneously converge in 
$C([-1,0])$ as $k\to\infty$. Furthermore, their limits $\vp_\infty:=\lim_k 
\vp_k$ and $\pi_i:=\lim_k p_i(\vp_k)$ have the property $\pi_i=p_i(\cdot , 
\vp_\infty)$, $i\in {\Bbb N}$. Thus, $\vp_\infty\in C_{a,\omega}^\infty 
([-1,0])$ and because of (\ref{2.3}), $\|\vp-\vp_\infty\|<2\ve$. 
\qed\medskip 

Let us recall that for any real number $a\neq 0$ there exists a unique 
solution $r$ to the equation  
\begin{eqnarray}\label{2.4} 
r(s)\equiv r(s;a)=1\chi_{[0,\infty)}(s)+a\int_0^s r(u-1)\chi_{[1,\infty)} 
(u)\, du\, , \quad s\ge 0. 
\end{eqnarray} 
In fact, we have 
\begin{eqnarray*} 
r(s)=0\ \ \mbox{\rm if}\ s<0\quad\mbox{\rm and}\quad r(s)=\sum_{l=0}^{k-1} 
\frac{a^l}{l!}(s-l)^l\ \ \mbox{\rm if}\ s\in [k-1,k),\ k\in {\Bbb N}. 
\end{eqnarray*} 
For $s< 0$, we set $r(s):=0$. Let us moreover collect some facts about 
the problem 
\begin{eqnarray}\label{2.5}
dX^{[-1,\infty)}_s=aX^{[-1,\infty)}_{s-1}+dW_s\, , \ s\ge 0,\quad X^{[-1, 
\infty)}_u=f(u),\ u\in [-1,0], 
\end{eqnarray} 
from \cite{Ri08}. For $s\ge 0$ and $u\in [-1,0]$ define $(I(s))(u)=0$ if 
$s+u<0$ and 
\begin{eqnarray*}
(I(s,W))(u)&&\hspace{-.5cm}:=W_{s+u}+\int_0^{s+u}W_{s-v+u}\, dr(v)-r(s+u) 
\cdot W_0 \\ 
&&\hspace{-.5cm}=\int_0^{s+u}r(s-v+u)\, dW_v\, , \quad s+u\ge 0. 
\end{eqnarray*} 
Also, there is a strongly continuous semigroup $(T_s)_{s\ge 0}$ in $C([-1,0 
])$ such that with $f\in C([-1,0])$, $T_sf(q)$ depends only on $f$ and $q+s 
$. In particular, for $s+q\in [-1,0]$, $T_sf(q)$ is directly defined by 
$T_sf(q):=f(s+q)$. For all other $s\ge 0$ and $q\in [-1,0]$, $T_sf(q)$ is 
implicitly given by the following. The function $g(v):=T_sf(v-s)$ for $v\in 
[-1,\infty)$ and any $s\in [v,v+1]$ is continuous on $v\in [-1,\infty)$, 
continuously differentiable on $v>0$, and satisfies  
\begin{eqnarray*}
g'(s)=ag(s-1)\, , \ s\ge 0,\quad g(u)=f(u)\, , \ u\in [-1,0]. 
\end{eqnarray*} 
There exists a pathwise unique solution $X^{[-1,\infty)}$ of (\ref{2.5}). 
In particular, 
\begin{eqnarray*}
X^{[-1,\infty)}_{s+u}\equiv X^{[-1,\infty)}_{s+u}(f,W)=(T_sf)(u)+(I(s,W)) 
(u)\, , \quad s\ge 0,\ u\in [-1,0]. 
\end{eqnarray*} 

Let $\omega$ be the restriction of $W\in\Omega$ to the interval $(-\infty 
,0]$. Below we will pathwise use the symbol $C_{a,W}^\infty ([-1,0])$ for 
$C_{a,\omega}^\infty([-1,0])$. 
\begin{lemma}\label{Lemma2.3}
Let $a\in (-1,0)$ and $v_0\in [0,1]$. (a) For $Q_{\sbnu}$-a.e. $W$ there 
exists a unique $f\equiv f(W,\cdot)\in C_{a,W}^\infty([-1,0])$ satisfying 
the following. 
\begin{itemize}
\item[(i)] For all $k\in {\Bbb N}$, the function $f$, considered as an 
element belonging to $C_{a,W}^k([-1,0])$, has the representation 
\begin{eqnarray*} 
f(v_0)&&\hspace{-.5cm}=\int_{-1}^{v_0} a\left(\ldots a\left(\int_{-1}^{v_{ 
k-1}}a\left(p_k(v_{k})-W_{-k-1}+W_{v_k-k}\vphantom{\dot{f}}\right)\, dv_k- 
W _{-k}\right.\right. \\ 
&&\hspace{-1.2cm}\left.\left.+W_{v_{k-1}-(k-1)}\vphantom{\int_0^{v_{k-1}}} 
\right)\, dv_{k-1}\ldots-W_{-2}+W_{v_1-1}\right)\, dv_1-W_{-1}+W_{v_0}    
\end{eqnarray*} 
for some sequence $p_k$ uniformly bounded in the $\sup$-norm with respect 
to $k\in {\Bbb Z}_+$ with $p_k(\cdot)-W_{-k-1}+W_{\cdot-k}\in C_{a,W_{\cdot 
-k}}^\infty([-1,0])$. 
\item[(ii)] The function $f$ has the representation 
\begin{eqnarray*} 
f(v_0)&&\hspace{-.5cm}=\sum_{k=1}^\infty a^k\int_{-1}^{v_0}\ldots\int_{-1 
}^{v_{k-1}}\left(-W_{-k-1}+W_{v_k-k}\vphantom{l^1}\right)\, dv_k\ldots dv_1 
 \\ 
&&\hspace{-.0cm}\vphantom{\sum}-W_{-1}+W_{v_0}\, ,
\end{eqnarray*} 
where the infinite sum converges absolutely. 
\end{itemize} 
\end{lemma} 
Proof.  We use 
\begin{eqnarray}\label{2.6}
\limsup_{k\to\infty}\sup_{-1\le v\le 0}\frac{-W_{-k-1}+W_{v-k}} 
{\D\sqrt{2\log k}}=1\quad Q_{\sbnu}\mbox{\rm -a.e.} 
\end{eqnarray} 
cf. \cite{MPSW10}, Exercise 5.1, from which it follows that the expression 
in (ii) is $Q_{\sbnu}$\mbox{\rm -a.e.} well-defined and absolutely converging 
for every $v_0\in [-1,0]$. From (ii) we get (i). 
\qed\medskip 
\begin{proposition}\label{Proposition2.4} 
Let $a\in (-1,0)$ and let $p_0,p_1,\ldots $ be the functions given by Lemma 
\ref{Lemma2.3} and $p_0:=f-(-W_{-1}+W)$. \\ 
(a) The infinite sum 
\begin{eqnarray*} 
q(v)\equiv q(W,v):=\sum_{k=0}^\infty\left(p_k(0)-W_{-k-1}+W_{-k}\vphantom 
{l^1}\right)\cdot r(v+k) 
\end{eqnarray*} 
converges absolutely $Q_{\sbnu}$-a.e. for any $v\in (-\infty,0]$. The random 
function $q:(-\infty,0]\to {\Bbb R}$ is $Q_{\sbnu}$-a.e. bounded on $(-\infty, 
0]$. \\ 
(b) The random function $X^{(-\infty,0)}\equiv X^{(-\infty,0)}(W):(-\infty 
,0)\to {\Bbb R}$ given by 
\begin{eqnarray*} 
X_v^{(-\infty,0)}=\sum_{k=0}^\infty\left(p_k(v+k)-W_{-k-1}+W_v\vphantom{l^1} 
\right)\cdot\chi_{[-k-1,-k)}(v)+q(v)\, ,  
\end{eqnarray*} 
$v\in (-\infty,0)$, is the unique random function that satisfies $Q_{\sbnu}$-a.e. 
\begin{itemize} 
\item[(i)] $dX_{s}=aX_{s-1}+dW_s$, $s\in (-\infty,0)$, and 
\item[(ii)] $\limsup_{v\to -\infty}|X_v^{(-\infty,0)}|<\infty$. 
\end{itemize} 
Moreover, we have 
\begin{itemize} 
\item[(iii)] $X^{(-\infty,0)}(W)=X^{(-\infty,0)}(W+x\1)$, $x\in {\Bbb R}$. 
\end{itemize} 
\end{proposition} 
Proof.  Part (a) follows from (\ref{2.6}) and that by $a\in (-1,0)$ 
the solution $r(u)$ of (\ref{2.4}) tends to $0$ as $u\to\infty$ exponentially 
fast. For this, recall that all roots of the characteristic equation 
$\lambda=ae^{-\lambda}$ have a negative real part bounded away from zero. 
For a detailed description of the behavior of $r(u)$ as $u\to\infty$ we 
refer to \cite{Er09}, \cite{GH01}, and \cite{HV93}. 

For part (b) we first check (i) on the open intervals $(-k-1,-k)$ stepwise 
relative to $k\in {\Bbb Z}_+$. To verify (the integrated version of) (i) at 
$-k$, $k\in {\Bbb N}$ we use (\ref{2.4}). Property (ii) of part (b) is a 
consequence of the $Q_{\sbnu}$-a.e. uniform boundedness of the $p_k$, $k\in 
{\Bbb Z}_+$, in the sup-norm, cf. Lemma \ref{Lemma2.3} (a), relation 
(\ref{2.6}), and the $Q_{\sbnu}$-a.e. boundedness of $q$, cf. part (a) of this 
proposition. Property (iii) is an immediate consequence of the definition 
of $X^{(-\infty,0)}$. 
\qed\medskip  
\medskip

By definition there is a continuous extension of $X^{(-\infty,0)}$ to a 
random function $(-\infty,0]\to {\Bbb R}$ which we denote by $X^{(-\infty, 
0]}$. 
\begin{theorem}\label{Theorem2.5} 
Let $a\in (-1,0)$. For $Q_{\sbnu}$-a.e. $W\in\Omega$ there exists a pathwise 
unique solution to the equation 
\begin{eqnarray*}
d\t X_s=a\t X_{s-1}\, ds+dW_s\, ,\quad s\in {\Bbb R},
\end{eqnarray*} 
satisfying $\limsup_{v\to -\infty}|\t X_v|<\infty$. The following 
holds. 
\begin{itemize} 
\item[(i)] The process $\t X$ is the solution to (\ref{2.1}) where 
\begin{eqnarray*}
\vp\equiv\vp (\, \cdot,W):=\left. X^{(-\infty,0]}(W)\right|_{[-1,0]}  
\end{eqnarray*} 
and $X^{(-\infty,0]}$ is the process of Proposition \ref{Proposition2.4} (b). 
\item[(ii)] $\t X_s=X_s^{(-\infty,0]}(W)$, $s\in (-\infty,0]$, and $\t X_s=X^{[-1, 
\infty)}_s(\vp,W)$, $s\in [-1,\infty)$. 
\item[(iii)] $\t X_{\cdot+v}(W)=\t X(W_{\cdot +v})$, $v\in {\Bbb R}$, and $\t X(W) 
=\t X(W+x\1)$, $x\in {\Bbb R}$.  
\end{itemize} 
\end{theorem}
Proof.  We first recall that any non-zero solution of $dx_s=ax_{s-1}\, 
ds$, $s\in {\Bbb R}$, satisfies $\limsup_{v\to -\infty}|x_v|=+\infty$. This is 
also a consequence of the fact that all roots of the characteristic equation 
$\lambda=ae^{-\lambda}$ have a negative real part bounded away from zero. Now 
the main statement of the theorem along with the properties (i) as well as 
(ii) follow from Proposition \ref{Proposition2.4} (b) and (\ref{2.5}). Here 
we also use the $Q_{\sbnu}$-a.e. uniqueness of the solution to (\ref{2.5}). 
The first part of property (iii) is obtained from the fact that uniqueness of 
$\t X$ is determined by $\limsup_{v\to -\infty}|\t X_v|<\infty$, the second 
part is obvious by the stochastic equation. 
\qed\medskip  
\medskip 

By time reversal we obtain the corresponding assertions relative to equation 
(\ref{2.2}). 
\begin{corollary}\label{Corollary2.6} 
Let $a\in (-1,0)$. For $Q_{\sbnu}$-a.e. $W\in\Omega$ there exists a pathwise 
unique solution to the equation 
\begin{eqnarray*}
d\t X_s=-a\t X_{s+1}\, ds+dW_s\, ,\quad s\in {\Bbb R},
\end{eqnarray*} 
satisfying $\limsup_{v\to\infty}|\t X^{(-\infty,0)}(v)|<\infty$. The following 
holds. 
\begin{itemize} 
\item[(i)] The process $\t X$ is the solution to (\ref{2.2}) where 
\begin{eqnarray*}
\vp\equiv\vp(u,W):=X_{-u}^{(-\infty,0]}(W_{-\,\cdot})\, ,\quad u\in [0,1], 
\end{eqnarray*} 
and $X^{(-\infty,0]}$ is the process of Proposition \ref{Proposition2.4} (b). 
\item[(ii)] $\t X$ is the time reversal of the stochastic process constructed in 
Theorem \ref{Theorem2.5} for $W_{-\, \cdot}$ instead of $W$. It satisfies (iii) 
of Theorem \ref{Theorem2.5}. 
\end{itemize} 
\end{corollary}

\section{Absolute Continuity under Time Shift} 
\setcounter{equation}{0}

\begin{theorem}\label{Theorem3.1}
Let 
\begin{eqnarray*}
b_0\equiv b_0(W):=f(0)+q(0)-W_0
\end{eqnarray*} 
where $f$ and $q$ are defined in Lemma \ref{Lemma2.3} and Proposition 
\ref{Proposition2.4}. Then $Q_{\sbnu}$-a.e. there exists a pathwise unique 
solution to the equation 
\begin{eqnarray}\label{3.1}
dX_s=a\left(X_{s-1}+b_0\vphantom{l^1}\right)\, ds+dW_s\, ,\quad s\in {\Bbb 
R},
\end{eqnarray} 
such that $\limsup_{v\to -\infty}|X_v|<\infty$. This solution 
is temporally homogeneous, i. e. for $A\equiv A(W):=X(W)-W$ and $W^u:=W_{\cdot 
+u}+A_u(W)\1$, $u\in {\Bbb R}$, it satisfies $A_0(W)=0$ and 
\begin{eqnarray*} 
X_{\cdot+v}(W)=X(W^v)\, , \quad v\in {\Bbb R}. 
\end{eqnarray*} 

\end{theorem}
Proof.  Let $\t X$ be the process introduced in Theorem \ref{Theorem2.5}. 
By $\limsup_{v\to -\infty}|\t X_v|<\infty$ there, the process 
\begin{eqnarray}\label{3.2}
X=\t X-b_0\1 
\end{eqnarray} 
is the unique solution to (\ref{3.1}) such that $\limsup_{v\to -\infty} 
|X_v|<\infty$. The definitions of $b_0$ and $A$ yield $b_0(W)=\t X_0(W)- 
W_0$ and $A_0(W)=0$ for $Q_{\sbnu}$-a.e. $W\in\Omega$. The definition of 
$b_0$ and Theorem \ref{Theorem2.5} (iii) imply furthermore 
\begin{eqnarray*}
b_0(W)&&\hspace{-.5cm}=\left(W_v+A_v(W)+b_0(0)\vphantom{l^1}\right)- 
\left(W_v+A_v(W)\vphantom{l^1}\right) \\ 
&&\hspace{-.5cm}=X_v(W)+b_0(W)-W_0^v=\t X_v(W)-W_0^v\vphantom{l^1} \\ 
&&\hspace{-.5cm}=\t X_0\left(W_{\, \cdot +v}\vphantom{l^1}\right)- 
W_0^v=\t X_0\left(W_{\, \cdot +v}+A_v(W)\1\vphantom{l^1}\right)-W_0^v  
 \\ 
&&\hspace{-.5cm}=\t X_0(W^v)-W_0^v =b_0(W^v)\, ,\quad v\in {\Bbb R}. 
\vphantom{l^1}
\end{eqnarray*}
Now we get from Theorem \ref{Theorem2.5} (iii) 
\begin{eqnarray*}
X_{\, \cdot +v}(W)&&\hspace{-.5cm}=\t X_{\, \cdot +v}(W)-b_0(W)\1 
\vphantom{l^1} \\ 
&&\hspace{-.5cm}=\t X\left(W_{\, \cdot +v}\right)-b_0(W)\1\vphantom 
{l^1} \\ 
&&\hspace{-.5cm}=\t X\left(W_{\, \cdot +v}+A_v(W)\1\vphantom{l^1} 
\right)-b_0(W^v)\vphantom{l^1} \\ 
&&\hspace{-.5cm}=X(W^v)\, ,\quad v\in {\Bbb R}. \vphantom{l^1}
\end{eqnarray*}
\qed\medskip 
\medskip

An immediate consequence of time reversal is the following. 
\begin{corollary}\label{Corollary3.2} 
Let $a\in (-1,0)$ and let $b_0\equiv b_0(W)$ be the random variable defined 
in Theorem \ref{Theorem3.1}. Then $Q_{\sbnu}$-a.e. there exists a pathwise 
unique solution to the equation 
\begin{eqnarray*}
dX_s=-a\left(X_{s+1}+b_0\vphantom{l^1}\right)\, ds+dW_s\, ,\quad s\in {\Bbb 
R},
\end{eqnarray*} 
such that $\limsup_{v\to +\infty}|X_v|<\infty$. This solution is temporally 
homogeneous. 
\end{corollary}

We obtain the following quasi-invariance under time shift of trajectories. 
\begin{theorem}\label{Theorem3.3}
Let $m$ denote the density of $W_0$. Let $X$ be the Ornstein-Uhlenbeck type 
process established in Theorem \ref{Theorem3.1} or the process of Corollary 
\ref{Corollary3.2}. Denoting by $P_{\sbmu}$ the law of $X$ we have 
$P_{\sbmu}$-a.e. 

\begin{eqnarray*}
\rho_{-t}(X)=\frac{m(X_{-t})}{m(X_0)}\, ,\quad t\in {\Bbb R}. 
\end{eqnarray*} 
\end{theorem} 
Proof.  {\it Step 1 } Let us first focus on the case of the process 
with delay established in Theorem \ref{Theorem3.1}. Let the measure $Q_{x}$ be 
obtained from $Q_{\sbnu}$ by conditioning on $W_0=x$, $x\in {\Bbb R}$. We 
have $Q_{\sbnu}=\int Q_x\, m(x)\, dx$. Next we use property (iii) of 
Theorem \ref{Theorem2.5} and relation (\ref{3.2}). We have for $B=B_i\times 
B_p$ with $B_i\in {\mathcal B}({\Bbb R})$ and $B_p\in {\mathcal F}\cap\{X\in 
\Omega:X_0=0\}$, 
\begin{eqnarray}\label{3.3}
&&\hspace{-.5cm}Q_{\sbnu}\left(W:X(W)\in B\right)=Q_{\sbnu}\left(W:\t X(W) 
-f(0)\1-q(0)\1+W_0\1\in B\right)\vphantom{l^1}\nonumber \\ 
&&\hspace{.5cm}=\int_{\{W_0\in B_i\}}Q_{W_0}\left(W:\t X(W)-f(0)\1-q(0)\1 
\in B_p\right)\, \bnu(dW_0) \vphantom{l^1}
\end{eqnarray} 
where we have used $X_0=W_0$ and $X=\t X-b_0\1=\t X-f(0)\1-q(0)\1+W_0\1$. 
Furthermore, it holds that $\t X(W)=\t X(W+x\1)$ according to Theorem 
\ref{Theorem2.5} (iii) and, by definition, $f(0)\equiv f(W,0)=f(W+x\1, 
0)$ as well as $q(0)\equiv q(W,0)=q(W+x\1,0)$, $x\in {\Bbb R}$. This says 
that $Q_x(W:\t X(W)-f(0)\1-q(0)\1\in B_p)$ is independent of $x$. For the 
next chain of equations we recall that $X_{\, \cdot -t}(W)=X\left(W^{-t} 
\right)$, cf. Theorem \ref{Theorem3.1}, $t\in {\Bbb R}$. We obtain 
\begin{eqnarray*}%\label{3.4}
&&\hspace{-.5cm}Q_{\sbnu}\left(W:X_{\, \cdot -t}(W)\in B\vphantom{l^1} 
\right)=Q_{\sbnu}\left(W:X\left(W^{-t}\right)\in B\vphantom{l^1}\right) 
\nonumber \\ 
&&\hspace{.0cm}=Q_{\sbnu}\left(W:\t X(W^{-t})-f(W^{-t},0)\1-q(W^{-t},0) 
\1+W^{-t}_0\1\in B\vphantom{l^1}\right)\nonumber \\ 
&&\hspace{.0cm}=Q_{\sbnu}\left(W:\t X(W^{-t})-f(W^{-t},0)\1-q(W^{-t},0) 
\1\in B_p,\, W^{-t}_0\1\in B_i\vphantom{l^1}\right) 
\end{eqnarray*} 
% \nonumber \\ 
% 
\begin{eqnarray}\label{3.4}
&&\hspace{.0cm}=\int_{x\in B_i}Q_{\sbnu}\left(W:\t X(W^{-t})-f(W^{-t},0) 
\1-q(W^{-t},0)\1\in B_p\left|W^{-t}_0=x\vphantom{l^1}\right.\right) 
\nonumber \\
&&\hspace{6.5cm}\, Q_{\sbnu}\left(W:W^{-t}_0\in dx\right)\nonumber \\ 
&&\hspace{.0cm}=\int_{x\in B_i}Q_{\sbnu}\left(W:\t X(W_{\, \cdot 
-t})-f(W_{\, \cdot -t},0)\1-q(W_{\, \cdot -t},0)\1\in B_p\left|W^{ 
-t}_0=x\vphantom{l^1}\right.\right)\nonumber \\
&&\hspace{6.5cm}\, Q_{\sbnu}\left(W:W^{-t}_0\in dx\right)\nonumber \\ 
&&\hspace{.0cm}=\int_{x\in B_i}Q_x\left(W:\t X(W)-f(W,0)\1-q(W,0)\1\in 
B_p\right)\, Q_{\sbnu}\left(W:W^{-t}_0\in dx\right) \qquad %\nonumber \\ 
\end{eqnarray} 
$t\in {\Bbb R}$. Recalling (\ref{1.2}), the claim follows now from 
(\ref{3.3}) and (\ref{3.4}) where use we 
that $Q_x(W:\t X(W)-f(0)\1-q(0)\1\in B_p\vphantom{l^1})$ is independent 
of $x$. We note also that $W_0=X_0$, $W^{-t}_0=X_{-t}(W)$, and that the 
distribution of $W_0$ is $\bnu$. 
\medskip

\nid
{\it Step 2 } Again, an immediate consequence of time reversal is the claim 
for he process with anticipation of Corollary \ref{Corollary3.2}. 
\qed\medskip 

\nid 
{\bf Remarks. } (1) It is also reasonable to verify the conditions of 
\cite{Lo13}, Theorem 1.12 and Proposition 1.16. This would lead to (\ref{1.3}), 
\begin{eqnarray*}
\rho_{-t}(X)&&\hspace{-.5cm}=\frac{\D m\left(X_{-t}\vphantom{l^1}\right)}{\D 
m\left(X_0\vphantom{l^1}\right)}\cdot\vphantom{\left\{\sum_{\tau_k}\right\}} 
\left|\nabla_{X_0}X_{-t}\vphantom{\dot{f}}\right| \\ 
&&\hspace{-.5cm}=\frac{\D m\left(X_{-t}\vphantom{l^1}\right)}{\D 
m\left(X_0\vphantom{l^1}\right)}\cdot\vphantom{\left\{\sum_{\tau_k}\right\}} 
\left|\nabla_{W_0}\left(\t X-f(0)\1-q(0)\1+W_0\1\vphantom{\dot{f}}\right) 
\right| \\ 
&&\hspace{-.5cm}=\frac{\D m\left(X_{-t}\vphantom{l^1}\right)}{\D m\left(X_0 
\vphantom{l^1}\right)}\, ,\quad t\in {\Bbb R}, 
\end{eqnarray*} 
where again we have used $\t X(W)=\t X(W+x\1)$, $f(W,0)=f(W+x\1,0)$, and 
$q(W,0)=q(W+x\1,0)$, $x\in {\Bbb R}$. The above proof demonstrates that, in 
comparably simple situations, the result can be obtained independently by 
direct and short calculations. 
\medskip 

\nid
(2) We observe the same simple form of the Radon-Nikodym derivative as in 
(\ref{1.1}).

%\par\bigskip\noindent
%{\bf Acknowledgment.} ????? {\bf(optional)}

\bibliographystyle{amsplain}

\end{document}